\documentclass[notitlepage,11pt,draft]{article}
\usepackage{amssymb,amsfonts,amsmath,geometry,esint,upgreek,cases,bbm,ulem}
\usepackage{scalerel,stackengine,mathtools,stmaryrd,soul,yhmath}

\usepackage{color,geometry}

\catcode`\@=11
\@addtoreset{equation}{section}

\catcode`\@=12

\allowdisplaybreaks

\def\proof{\noindent{\textbf{Proof. }}}
\def\QED{\hfill {$\square$}\goodbreak \medskip}

\newtheorem{Theorem}{Theorem}[section]
\newtheorem{Lemma}[Theorem]{Lemma}
\newtheorem{Remark}[Theorem]{Remark}
\newtheorem{Example}[Theorem]{Example}

\def\R{\mathbb{R}}
\def\f{\varphi}
\def\eps{\varepsilon}
\def\S{\mathbb S}
\def\disk{{\mathbb D^2}}
\def\ball{{\mathbb B^4}}

\def\rme{{\rm e}}
\def\Crme{{\mathcal C^\rme}}

\begin{document}

\title{\Large{Rellich-Poincar\'e inequalities and improved Rellich inequalities}}

\author{Paolo Caldiroli\footnote{Dipartimento di Matematica ``G. Peano", Universit\`a di Torino, Via Carlo Alberto 10, 10123 Torino, Italy. Email: {paolo.caldiroli@unito.it}. } ~and~
 Roberta Musina\footnote{Dipartimento di Scienze Matematiche, Informatiche e Fisiche, Universit\`a di Udine,
via delle Scienze 206, 33100 Udine, Italy. Email: {roberta.musina@uniud.it}}}

\date{}

\maketitle

\begin{abstract}
{We study linear variational problems on bounded domains $\Omega \subset \R^d$
involving the Dirichlet biharmonic operator $\Delta^{\!2}$ and the Rellich potential $|x|^{-4}$,
under the assumption $0\in\overline\Omega$.
The case $0\in\partial\Omega$ 
is particularly noteworthy and subtle.
{\small }}

\vskip0.5cm

\noindent
\textbf{Keywords:} {Rellich inequality, best constant, singular linear problems.}

\medskip\noindent
\textbf{2010 Mathematics Subject Classfication:} 26D10, 39B62, 45R05.
\end{abstract}

\date{}

\maketitle

\section{Introduction}
We are interested in linear variational problems on domains $\Omega \subset \R^d$ 
governed by the Dirichlet (clamped) biharmonic operator $\Delta^{\!2}$ and involving the Rellich potential $|x|^{-4}$.
We restrict our attention to bounded domains; the case of unbounded domains 
is studied in \cite{CM2026}.

To avoid trivial situations, we will always assume that either $0\in\Omega$ or $0\in\partial\Omega$, 
the latter being a particularly interesting and challenging case.

Our work was initially inspired by the celebrated result of Brezis and V\'azquez  \cite[Theorem 4.1]{BV},
where it is proved that 
for $0\in\Omega$ and $d\ge 3$ the following inequality holds: 
\begin{equation}
\label{eq:BV_original}
\int_\Omega|\nabla u|^2\,dx-\left(\frac{d-2}{2}\right)^2\int_\Omega\frac{|u|^2}{|x|^2}\,dx\ge \left(\frac{|B_1|}{|\Omega|}\right)^{\!\!\frac2d}
\lambda(-\Delta,\disk)\int_\Omega|u|^2\,dx~,\quad u\in H^1_0(\Omega)\,.
\end{equation}
Here $|\Omega|$, 
$|B_1|$ denote the volumes of $\Omega$ and of the unit ball $B_1\subset \R^d$, respectively, and
$\lambda(-\Delta,\disk)$
is the Poincar\'e constant on the unit disk $\disk\subset \R^2$, namely
\begin{equation}
\label{eq:disk_lambda}
\lambda(-\Delta,\disk)=\inf_{v\in H^1_0(\disk)\atop v\neq 0}
\frac{\displaystyle\int_{\disk}|\nabla v|^2\,dx}{\displaystyle\int_{\disk}| v|^2\,dx}\,.
\end{equation}
We  also refer to \cite{FM, Fall}, where the case $0\in\partial\Omega$ is studied, and to \cite{Fall2} for a related problem. 

\medskip

Let us describe the problem we are interested in and our main results (see Section \ref{S:Preliminaries} for the notation and main definitions).
Throughout the paper, we assume that at least one of the following conditions is satisfied:
$$
(A1)~~d\ge 5~,\qquad (A2)~~0\in\partial \Omega~\text{and $\Omega$ has an external direction at $0$.}
$$

This allows us to prove in Lemma \ref{L:mu0_Omega}, that  the {\it Rellich constant of $\Omega$}, which is defined by
\begin{equation}
\label{eq:mu0_Omega}
\mu_0(\Omega)= \inf_{u\in H^2_0(\Omega)\atop u\neq 0}
\frac{\displaystyle\int_{\Omega}|\Delta u|^2\,dx}{\displaystyle\int_{\Omega}|x|^{-4}| u|^2\,dx}\,,
\end{equation}
is positive.

If $d\ge 5$, then the weight $|x|^{-4}$ is locally integrable on $\R^d$, and the Rellich inequality \cite{Re} yields
\begin{equation}
\label{eq:Rellich_space}
\mu_0(\R^d)= \inf_{u\in C^\infty_c(\R^d)\atop u\neq 0}
\frac{\displaystyle\int_{\R^d}|\Delta u|^2\,dx}{\displaystyle\int_{\R^d}|x|^{-4}| u|^2\,dx}= \gamma_d^2~,
\qquad \text{where}\quad \gamma_d=\frac{d(d-4)}{4}\,.
\end{equation}
It follows that $\mu_0(\Omega)\ge \gamma_d^2>0$. 
If $0\in\Omega$, then
$\mu_0(\Omega)=\mu_0(\R^d)=\gamma_d^2$ and 
$\mu_0(\Omega)$ is not achieved in $H^2_0(\Omega)$, because 
of the scale invariance of the Rayleigh quotient in (\ref{eq:Rellich_space}).

If $(A2)$ is satisfied, then no restriction on the dimension $d$ is needed,
and in fact $\mu_0(\Omega)>0$ even in low dimensions. 
In this case, $\mu_0(\Omega)$ may be attained, see the examples in Appendix \ref{A:Example}.

\medskip

We  then introduce, following \cite{BM, FM}, the minimization problem
\begin{equation}
\label{eq:mu_def}
\mu_\lambda(\Omega)= \inf_{u\in H^2_0(\Omega)\atop u\neq 0} J_\lambda(u)~,
\qquad 
J_\lambda(u):=\frac
{\displaystyle\int_\Omega|\Delta u|^2~\!dx-\lambda\displaystyle\int_\Omega| u|^2~\!dx}
{\displaystyle\int_\Omega|x|^{-4}| u|^2~\!dx}\,.
\end{equation}

\medskip

If $\lambda$ coincides with the principal eigenvalue of  the  operator $\Delta^{\!2}$ on $H^2_0(\Omega)$, then
 $\mu_\lambda(\Omega)=0$ and  the infimum $\mu_\lambda(\Omega)$ is achieved by any correspondent eigenfunction. We can therefore define
$$
\lambda^*_\Omega=\inf\big\{\lambda\in\R~|~\mu_\lambda(\Omega)~~\text{is achieved}\big\}\in[-\infty,\infty)\,.
$$
In our first result we describe the main features 
of the function $\lambda\mapsto \mu_\lambda(\Omega)$ and of the minimization problem in (\ref{eq:mu_def}). 
The variational problem defining $\mu_\lambda(\Omega)$ is not compact because of the action of dilations in $\R^d$, and in fact a minimizing sequence for $\mu_\lambda(\Omega)$ 
may concentrate at the origin.
The critical threshold is given by
$$
\mu^\infty_\Omega=\lim_{\lambda\to \lambda^*_\Omega}\mu_\lambda(\Omega)\,.
$$
Indeed, in our first result we show, in particular, that 
$\mu_\lambda(\Omega)\le \mu^\infty_\Omega$ for any
$\lambda\in\R$, and that the validity of the strict inequality $\mu_\lambda(\Omega)< \mu^\infty_\Omega$ guarantees the existence of a minimizer
for $\mu_\lambda(\Omega)$.

\begin{Theorem} 
\label{T:in_bordo}
Let $\Omega\subset\R^d$, $d\ge 1$, be a bounded domain. Assume that $(A1)$ or $(A2)$ is satisfied. Then
$\lambda^*_\Omega<\lambda(\Delta^{\!2},\Omega)$, 
$\mu^\infty_\Omega>0$, and 
the following facts hold:
\begin{itemize}
\item[$i)$] If $\lambda>\lambda^*_\Omega$, then $\mu_\lambda(\Omega)<\mu^\infty_\Omega$ and 
$\mu_\lambda(\Omega)$ is achieved;
\item[$ii)$] The function $\lambda\mapsto \mu_\lambda(\Omega)$ is strictly decreasing and strictly concave in $(\lambda^*_\Omega,\infty)$;
\item[$iii)$] If, in addition, $\lambda^*_\Omega>-\infty$, then
$\mu_\lambda(\Omega)=\mu^\infty_\Omega$ for any $\lambda\le \lambda^*_\Omega$. 
\end{itemize}
\end{Theorem}
In Appendix \ref{A:Example} we show that there exist domains $\Omega$ satisfying $(A2)$ and such that 
$\mu^\infty_\Omega=\infty$, which implies that $\lambda^*_\Omega=-\infty$.

Further information can be obtained either in case of an interior singularity 
or when $\Omega$ is sufficiently smooth at $0\in\partial\Omega$. 
To state our main results, it is convenient to distinguish between these two cases.

If $0\in\Omega$, which needs $d\ge 5$,  we can provide a positive lower bound for $\lambda^*_\Omega$.
This readily implies an
improved Rellich-type inequality in the spirit of \cite{BV}. Our bound involves
the radius
$$
R=\sup_{x\in\Omega}|x|~,
$$
the Poincar\'e constant on the unit disk already appeared in (\ref{eq:BV_original}), 
and the principal eigenvalue of the biharmonic operator $\Delta^{\!2}$ on the unit $4$-dimensional ball with clamped boundary conditions,
namely, 
\begin{equation}
\label{eq:ball_lambda}
\lambda(\Delta^{\!2},\ball)=\inf_{v\in H^2_0(\ball)\atop v\neq 0}
\frac{\displaystyle\int_{\ball}|\Delta v|^2\,dx}{\displaystyle\int_{\ball}| v|^2\,dx}\,.
\end{equation}

\begin{Theorem}
\label{T:in_Omega}
Let $\Omega$ be a bounded 
domain in $\R^d$, $d\ge 5$, such that $0\in \Omega$. 
Then $\mu^\infty_\Omega=\mu_0(\R^d)$ and 
\begin{equation}
\label{eq:lambda*1}
\lambda^*_\Omega\ge \frac{\lambda(\Delta^{\!2},\ball)}{R^4}+2d(d-4)\frac{\lambda(-\Delta,\disk)}{R^4}\,.
\end{equation}
Therefore,  for any $u\in H^2_0(\Omega)$ 
the following inequality holds:
\begin{equation}
\label{eq:impr1}
\int_{\Omega}|\Delta u|^2~\!dx- \left(\frac{d(d-4)}4\right)^2
\int_{\Omega}\frac{|u|^2}{|x|^4}~\!dx\ge 
\left[\frac{\lambda(\Delta^{\!2},\ball)}{R^4}+2d(d-4)\frac{\lambda(-\Delta,\disk)}{R^4}\right]\int_{\Omega}| u|^2\,dx\,.
\end{equation}
\end{Theorem}
We conjecture that the infimum $\mu_\lambda(\Omega)$ is achieved if and only if $\lambda>\lambda^*_\Omega$.
Proving nonexistence at the critical value $\lambda=\lambda^*_\Omega$ is quite challenging in this setting,
since minimizers of $\mu_\lambda(\Omega)$ do not have, in general, a constant sign.

Inequality (\ref{eq:impr1}) is complementary to the inequality of Gazzola--Grunau--Mitidieri \cite{GGM}, namely, 
$$
\int_{\Omega}|\Delta u|^2~\!dx-\left(\frac{d(d-4)}4\right)^2
 \int_{\Omega}\frac{|u|^2}{|x|^4}~\!dx\ge 
\frac{\lambda(\Delta^{\!2},\ball)}{R^4}\int_{\Omega}| u|^2~\!dx
+ \frac{d(d-4)}{2}\frac{\lambda(-\Delta,\disk)}{R^2}\int_{\Omega}\frac{|u|^2}{|x|^2}~\!dx
$$
although the two estimates are not directly comparable.

\medskip

In case of a boundary singularity, the critical threshold $\mu^\infty_\Omega$ depends on the geometry of
$\partial \Omega$ at $0$. It turns out that for sufficiently smooth domains, $\mu^\infty_\Omega$ coincides with 
$$
\mu_0(\R^d_+)=
 \inf_{u\in C^\infty_c(\R^d_+)\atop u\neq 0}
\frac{\displaystyle\int_{\R^d_+}|\Delta u|^2\,dx}{\displaystyle\int_{\R^d_+}|x|^{-4}| u|^2\,dx}\,.
$$
By the results in \cite{CM}, we have that $\mu_0(\R^d_+)$ is positive in any dimension, although for $d\ge 2$, its explicit value is not known. 

For the sake of clarity, we state our results separately for the cases $d\geq 2$ and $d=1$. 
The following result is the 
main contribution of the present paper.

\begin{Theorem} 
\label{T:smooth}
Let $\Omega$ be a bounded 
domain in $\R^d$, $d\ge 2$, such that $0\in \partial\Omega$. Assume that $\partial\Omega$ is differentiable at $0$.
Then  
\begin{equation}
\label{eq:limit}
\mu^\infty_\Omega=\mu_0(\R^d_+)\,.
\end{equation}
If in addition $\Omega$ is contained in a half-space, then
\begin{equation}
\label{eq:estimate_half}
\lambda^*_\Omega\ge \dfrac{\lambda(\Delta^{\!2},\ball)}{R^4}+2(d^2-4)\dfrac{\lambda(-\Delta,\disk)}{R^4}\,.
\end{equation}
In particular, $\mu_0(\Omega)=\mu_0(\R^d_+)$, the Rellich constant of $\Omega$ is not achieved, and  for any $u\in H^2_0(\Omega)$ 
$$
\int_{\Omega}|\Delta u|^2~\!dx- 
\mu_0(\R^d_+) 
 \int_{\Omega}\frac{|u|^2}{|x|^4}~\!dx\ge 
\Big(\frac{\lambda(\Delta^{\!2},\ball)}{R^4}+2(d^2-4)\dfrac{\lambda(-\Delta,\disk)}{R^4}\Big)\int_{\Omega}| u|^2~\!dx\,.
$$
\end{Theorem}
By analogy with \cite{Fall}, we conjecture that $\lambda^*_\Omega$ is always finite, at least when $\Omega$ is sufficiently smooth.
We also conjecture that the infimum $\mu_{\lambda}(\Omega)$ is not attained  if $\lambda= \lambda^*_\Omega$, as was
proved in \cite{FM} for the corresponding lower-order problem. For both questions, the arguments in \cite{Fall}, \cite{FM} do not apply
in the present setting, due to the failure of the maximum principle, and to severe additional technical difficulties.
Finally, we notice that the improved Rellich inequality we obtain in case $\Omega\subset \R^d_+$
can be seen as the higher order version of \cite[Theorem 5.2]{FM}.

\medskip

The lowest dimensional case $d=1$ is exceptional, since the origin 
disconnects the real line. Up to reflection, it is  enough to study the minimization problem
\begin{equation}
\label{eq:mu_Uno}
\mu_\lambda(0,R)= \inf_{u\in H^2_0(0,R)\atop u\neq 0} \frac
{\displaystyle\int_0^R|u''|^2~\!dx-\lambda\displaystyle\int_0^R| u|^2~\!dx}
{\displaystyle\int_0^Rx^{-4}| u|^2~\!dx}\,.
\end{equation}
In our last result we show in particular that $\mu^\infty_{(0,R)}=\mu_0(\R_+)$. 
The Rellich constant $\mu_0(\R_+)$ can be explicitly computed: by applying the weighted Hardy inequality twice,
see (\ref{eq:Hardy1}),
one readily obtains 
\begin{equation}
\label{eq:916}
\mu_0(\R_+):=\inf_{u\in C^\infty_c(\R_+)\atop u\neq 0}
\frac{\displaystyle\int_0^\infty|u''|^2~\!dx}
{\displaystyle\int_0^\infty x^{-4}| u|^2~\!dx}=\frac{9}{16}\,.
\end{equation}
Our result for (\ref{eq:mu_Uno}) involves the eigenvalues in (\ref{eq:disk_lambda}), (\ref{eq:ball_lambda}), and   
the principal eigenvalue of the Laplace operator $-\Delta$ 
on $H^1_0(\ball)$,  which is denoted by $\lambda(-\Delta,\ball)$.

\begin{Theorem} 
\label{T:Uno}
Consider the domain $(0,R)\subset \R$. Then $\mu^\infty_{(0,R)}={9}/{16}$, and
\begin{equation}
\label{eq:estimate*_Uno}
\lambda^*_{(0,R)}\ge \dfrac{5}{8R^4}\,{\lambda(\Delta^{\!2},\ball)} +\dfrac{3}{4R^4}\,{\lambda(-\Delta,\ball)\lambda(-\Delta,\disk)}\,.
\end{equation} 
Therefore, for any $u\in H^2_0(0,R)$ the following inequality holds: 
$$
\int_0^R|u''|^2~\!dx-\frac{9}{16}\int_0^R\frac{|u|^2}{x^4}~\!dx\ge
\left(\dfrac{5}{8R^4} \lambda(\Delta^{\!2},\ball)+\dfrac{3}{8R^4}\lambda(-\Delta,\ball)\lambda(-\Delta,\disk)\right)
\int_0^R| u|^2~\!dx\,.
$$
\end{Theorem}

\medskip

For the convenience of the reader, we provide  here the approximate values of the 
constants appearing in our improved Rellich inequalities. 
It is well known that $\lambda(-\Delta,\disk)=j_{0,1}^2$ and $\lambda(-\Delta,\ball)=j_{1,1}^2$, where $j_{\nu,1}$ is the first positive
zero of the Bessel function $J_\nu$. For the approximate value of $\lambda(\Delta^{\!2},\ball)$ we cite \cite{Ta}.
We have
$$
\lambda(-\Delta,\disk)\approx 5.78\,,\quad \lambda(-\Delta,\ball)\approx 14.68\,,\quad \lambda(\Delta^{\!2},\ball)\approx 452.
$$

{The paper is organized as follows.

After the preliminary Section \ref{S:Preliminaries}, in Section \ref{S:general} we start to study the infimum $\mu_\lambda(\Omega)$.
In particular, we prove the compactness Lemma \ref{L:compactness}, which is a crucial step in our arguments.
Finally,  we prove Theorem \ref{T:in_bordo}.

In Section \ref{S:conic} we deal with the critical parameter  
$\lambda^*_\Omega$ in case $\Omega$ is a {\it conic domain}, namely, it is the intersection of a cone with a ball about the origin. 
In order to prove Lemma \ref{L:conic_domains1}, which is the main step in the proof of our improved Rellich inequalities,
we recall a few facts about the Emden-Fowler transform and study 
some auxiliary variational problems. The section concludes with the proofs of Theorem \ref{T:in_Omega}  and of
Theorem \ref{T:Uno}.

Section \ref{S:proofs} contains the proof of Theorem \ref{T:smooth}.

Appendix \ref{A:Example} contains two examples.
In Appendix \ref{A:BV} we prove two inequalities for radial functions 
in $H^2_0(\ball)$.
These inequalities are needed to prove our estimates on $\lambda^*_\Omega$.

\section{Notation and preliminaries}
\label{S:Preliminaries}

\begin{description}

\item$\bullet$ We denote by $B_R\subset \R^d$ the open Euclidean ball of radius $R>0$ about the origin. In the exceptional 
dimensions $d=2, 4$ we will sometimes use the notation $\disk$, $\ball$, respectively, instead of $B_1$.

\item$\bullet$ $\Omega$ is a domain in $\R^d$. If $0\in \partial\Omega$, we say that $\Omega$ has an external direction at $0$
if there exist $\delta>0$ such that $\Omega\cap B_\delta$ is contained in the complement of a half-line. That is,
if there exists
$\rme\in\S^{d-1}$ such that
$$
\Omega\cap B_\delta \subset \big\{r\sigma\in\R^d~|~r>0~,~~\sigma\in\S^{d-1}\setminus\{\rme\}\big\}\,.
$$
We say that {\it $\partial\Omega$ is differentiable at $0$}
if there exists ${\delta_0}>0$ and a  
continuous function $F\colon {B}_{\delta_0}\to\R$, differentiable at $0$, and such that
\begin{equation}
\label{eq:F_bordo}
|\nabla F(0)|=1~,\quad 
\Omega\cap B_{\delta_0}=\{x\in B_{\delta_0}\colon F(x)<0\}~, \quad \partial\Omega\cap B_{\delta_0}=\{x\in B_{\delta_0}\colon F(x)=0\}\,.
\end{equation}

\item$\bullet$ 
The norm in $L^2(\Omega)$ is denoted by $\|\cdot\|_2$. We use standard definitions and norms for weighted $L^2$-spaces and for the Sobolev spaces 
$H^1_0(\Omega)$, $H^2_0(\Omega)$.  
Functions in $H^1_0(\Omega)$ and   $H^2_0(\Omega)$ are extended by $0$ outside $\Omega$, so that 
these spaces are identified with subspaces of $H^1(\R^d)$ and  $H^2(\R^d)$, respectively.

Since
the $(m,2)$-capacity of the set $\{0\}$ vanishes if $2m\le d$, see for instance \cite[Theorem 9.9.1]{AH},
the next lemma hold.

\begin{Lemma}
\label{L:density}
If $d\ge 2$, then $C^\infty_c(B_1\setminus\{0\})$ is dense in $H^1_0(B_1)$.
If $d\ge 4$, then $C^\infty_c(B_1\setminus\{0\})$ is dense in $H^2_0(B_1)$.
\end{Lemma}

\item$\bullet$ if $d\ge 2$, we denote by 
$\nabla_{\!\sigma}, \Delta_\sigma$ and $d\sigma$ the gradient, the Laplace-Beltrami operator, and the integration measure on the unit sphere $\S^{d-1}$,
respectively.
For any domain $\omega\subseteq\S^{d-1}$, the spaces $H^1_0(\omega), H^2_0(\omega)$ are defined in  standard way.
We introduce the eigenvalue
\begin{equation}
\label{eq:E_omega}
{\rm E}(\omega)=\inf_{\f\in H^1_0(\omega)\atop \f\neq 0}
\frac{\displaystyle\int_\omega|\nabla_{\!\sigma}\f|^2\,d\sigma}{\displaystyle\int_\omega|\f|^2\,d\sigma}\,.
\end{equation}
The infimum ${\rm E}(\omega)$ is achieved by Rellich's theorem. 
For instance, ${\rm E}(\S^{d-1})=0$ as
 $H^1_0(\S^{d-1})=H^1(\S^{d-1})$; 
if $\S^{d-1}_+$ is any hemisphere, then ${\rm E}(\S^{d-1}_+)=d-1$.

\item$\bullet$ For a domain $\omega\subseteq \S^{d-1}$ we set
\begin{equation}
\label{eq:m_omega}
m_{\gamma_{d}}(\omega)= \inf_{\f\in H^2_0(\omega)\atop \f\neq 0}
\frac{\displaystyle\int_\omega|\Delta_\sigma\f-\gamma_{d}\f|^2\,d\sigma}{\displaystyle\int_\omega|\f|^2\,d\sigma}\,,
\end{equation}
where $\gamma_d=d(d-4)/4$, see (\ref{eq:Rellich_space}).
The infimum $m_{\gamma_{d}}(\omega)$ is achieved by  Rellich's  theorem. If $d\ge 5$, then evidently
$$
m_{\gamma_{d}}(\omega)\ge m_{\gamma_{d}}(\S^{d-1})=\gamma_d^2>0\,.
$$
Further, we have $m_{\gamma_{d}}(\S^{3})=0$,
$m_{\gamma_{d}}(\S^{2})=9/16$ and $m_{\gamma_{d}}(\S^{1})=0$, see \cite{CM}. 

\item$\bullet$ Let $\rme\in\S^{d-1}$. By Lemma \ref{L:density}, if $d\ge 5$ then $H^2_0(\S^{d-1}\setminus\{\rme\})=H^2(\S^{d-1})$, hence
$m_{\gamma_{d}}(\S^{d-1}\setminus\{\rme\})=m_{\gamma_{d}}(\S^{d-1})=\gamma_d^2$ in this case. 
If $d\in\{2,3,4\}$ we have that
$H^2_0(\S^{d-1}\setminus\{\rme\})\subsetneq H^2(\S^{d-1})$, hence
$$
m_{\gamma_{d}}(\S^{d-1}\setminus\{\rme\})>m_{\gamma_{d}}(\S^{d-1})\ge 0\,.
$$
Thus
in any dimension, $m_{\gamma_{d}}(\omega)>0$ whenever $\omega$ is  properly contained in $\S^{d-1}$. 

\item$\bullet$~A {\it cone} is an open, dilation-invariant domain $\mathcal C\subset \R^d$. The Rellich constant of the cone $\mathcal C$ is 
$$
\mu_0(\mathcal C):=\inf_{u\in C^\infty_c(\mathcal C)\atop u\neq 0}
\frac{\displaystyle\int_\mathcal C|\Delta u|^2\,dx}{\displaystyle\int_\mathcal C|x|^{-4}| u|^2\,dx}\,.
$$
If $d\ge 2$, we set
\begin{equation}
\label{eq:omegaC}
\omega_\mathcal C=\mathcal C\cap S^{d-1}\,,
\end{equation}
and we say that $\mathcal C$ is a {\it proper cone} if $\mathcal C\subsetneq\R^d\setminus\{0\}$.
In this case, 
\begin{equation}
\label{eq:mu_m}
\mu_0(\mathcal C)= m_{\gamma_{d}}(\omega_\mathcal C)>0\,.
\end{equation}
The equality in (\ref{eq:mu_m}) can be proved by adapting the proof of \cite[Theorem 3.1]{CM}. In particular, denoting by $\R^d_+$ any half space in $\R^d$,  
we have that $\mu_0(\R^d_+)>0$. 

The only proper cones in $\R$ are the half lines $\R_\pm=\{\pm x>0\}$. Recall that for any $\alpha\in\R$, the weighted Hardy inequality
\begin{equation}
\label{eq:Hardy1}
\int_0^\infty x^\alpha |f'|^2\,dx\ge \Big(\frac{\alpha-1}{2}\Big)^2\int_0^\infty x^{\alpha-2}|f|^2\,dx~,\qquad f\in C^\infty_c(\R_+)\,,
\end{equation}
holds with a sharp constant in the right hand side.

\item$\bullet$ $c>0$ denotes a generic constant independent of the relevant quantities, whose value may change from line to line.

\end{description}

\section{A first look on $\mu_\lambda(\Omega)$ and proof of Theorem \ref{T:in_bordo}}
\label{S:general}

We start by studying, for any domain $\Omega\subset \R^d$, the {\it Rellich constant} $\mu_0(\Omega)$ in (\ref{eq:mu0_Omega}) 
and the {\it Hardy-Rellich constant}
$$
\mu_0^{\rm HR}(\Omega):=\inf_{u\in H^2_0(\Omega)\atop u\neq 0}
\frac{\displaystyle\int_\Omega|\Delta u|^2dx}
{\displaystyle\int_\Omega|x|^{-2}|\nabla u|^2dx},
$$
compare with \cite{Be, BC, C, GM}. 

\begin{Lemma}
\label{L:mu0_Omega0}
If $\mu_0(\Omega)>0$, then $\mu_0^{\rm HR}(\Omega)>0$.
\end{Lemma}

\proof
For $u\in C^\infty_c(\Omega)$ we  estimate
$$
\begin{aligned}
\int_\Omega|x|^{-2}|\nabla u|^2\,dx&=-\int_\Omega|x|^{-2}u\Delta u\,dx-(d-4)\int_\Omega|x|^{-4}|u|^2\,dx\\
&\le \big(\int_\Omega|x|^{-4}| u|^2\,dx\Big)^\frac12\big(\int_\Omega|\Delta u|^2\,dx\Big)^\frac12
+\frac{(4-d)^+}{\mu_0(\Omega)}\int_\Omega|x|^{-2}|\Delta u|^2\,dx\\
&\le
\left(\frac{1}{\sqrt{\mu_0(\Omega)}}+\frac{(4-d)^+}{\mu_0(\Omega)} \right)\int_\Omega|\Delta u|^2\,dx\,.
\end{aligned}
$$
Since $C^\infty_c(\Omega)$ is dense in $H^2_0(\Omega)$, this ends the proof.
\QED

\begin{Lemma}
\label{L:mu0_Omega} 
 Let  $\Omega\subset\R^d$ be a bounded domain satisfying $(A1)$ or $(A2)$. The following facts hold.
\begin{itemize}
\item[$i)$] $\mu_0(\Omega)>0$;
\item[$ii)$] $\mu_\lambda(\Omega)$ is finite, for any $\lambda\in\R$;
\item[$iii)$] The function $\lambda\mapsto \mu_\lambda(\Omega)$ is non increasing and concave. In particular, it is locally Lipschitz continuous.
\end{itemize}
\end{Lemma}

\proof
If $d\ge 5$ then $\mu_0(\Omega)\ge \mu_0(\R^d)=\gamma_d^2>0$.

If $d=1$ only the case $\Omega=(0,R)$ has to be considered.  Using dilations and the Hardy inequalities (\ref{eq:Hardy1}),
it is easy to show that $\mu_0(0,R)= \mu_0(\R_+)=9/16$.

\medskip

To complete the proof of $i)$ assume that $d\in\{2,3,4\}$ and that  $(A2)$ is satisfied. We can find ${\rm e}\in\S^{d-1}$ and
$\delta>0$, such that
$$\Omega\cap B_\delta\subset \mathcal C^{\rm e}\quad \text{where}\quad \Crme:=\big\{r\sigma\in\R^d~|~r>0~,~~\sigma\in\S^{d-1}\setminus\{\rme\}\big\}\,.
$$
Hence $\mu_0(\Omega\cap B_{2\delta})\ge \mu_0(\Crme)$
and $\mu_0(\Crme)>0$  by (\ref{eq:mu_m}).

We now take a cut-off function $\psi\in C^\infty_c(B_{2\delta})$ such that $0\le \psi\le 1$, $\psi\equiv 1$ in $B_\delta$. Fix $u\in C^\infty_c(\Omega)$.
Since $\psi u\in C^\infty_c(\Omega\cap B_{2\delta})$, using the continuity of the embeddings $H^2_0(\Omega)\hookrightarrow H^1_0(\Omega) \hookrightarrow L^2(\Omega)$, we can estimate
$$
\sqrt{\mu_0(\Omega\cap B_{2\delta})} \||x|^{-2}(\psi u)\|_2 \le \|\Delta(\psi u)\|_2 \le \|\Delta u\|_2+2\|\nabla \psi\cdot\nabla u\|_2+\|\Delta\psi\|_2\|u\|_2
\le c\|\Delta u\|_2,
$$
where $c>0$ depends only on $\Omega$ and $\psi$. Since $\text{supp}(1-\psi)\subseteq \{|x|\ge \delta\}$, we have
$$
\||x|^{-2}(1-\psi)u\|_2
\le \delta^{-2}\|u\|_2\le c\|\Delta u\|_2\,,
$$
for another constant $c>0$ depending only on $\delta$ and $\Omega$. 
Finally, since $u$ was arbitrarily chosen in $C^\infty_c(\Omega)$ and 
thanks to the triangle inequality
$\||x|^{-2}u\|_2\le \||x|^{-2}(\psi u)\|_2+\||x|^{-2}(1-\psi)u\|_2$,
we can conclude that $\mu_0(\Omega)>0$.

\medskip

Next, we take $R>0$ large, so that $\Omega\subseteq B_R$ and fix $\lambda\in\R$. For any $u\in H^2_0(\Omega)$ we have  
 $\|u\|_2\le R^{-2}\||x|^{-2}u\|_2$, 
 and therefore $J_\lambda(u)\ge  \mu_0(\Omega) -|\lambda|R^4$. Since $u$ was arbitrarity chosen, we infer that 
 $\mu_\lambda(\Omega)  \ge \mu_0(\Omega) -|\lambda|R^4>\infty$, as claimed in $ii)$.

Lastly, notice that for any $u\in H^2_0(\Omega)\setminus\{0\}$ fixed, the quantity $J_\lambda(u)$ is an affine and decreasing function of $\lambda$.
This readily gives  $iii)$.
\QED

The following lemma will play a crucial role in proving our main results.

\begin{Lemma}[compactness]
\label{L:compactness}
Let $\Omega\subset\R^d$, $d\ge 1$, be a bounded domain satisfying $(A1)$ or $(A2)$. 
Let $\lambda_0<\lambda_1$ be given, so that $\mu_{\lambda_1}(\Omega)\le \mu_{\lambda_0}(\Omega)$ by Lemma \ref{L:mu0_Omega}.
\begin{itemize}
\item[$i)$] If $\mu_{\lambda_1}(\Omega)< \mu_{\lambda_0}(\Omega)$, then $\mu_{\lambda}(\Omega)$ is achieved for any $\lambda\ge \lambda_1$;
\item[$ii)$] If $\mu_{\lambda_0}(\Omega)$ is achieved, then $\mu_{\lambda_1}(\Omega)< \mu_{\lambda_0}(\Omega)$, hence $\mu_{\lambda_1}(\Omega)$ is achieved;
\item[$iii)$] If $\mu_{\lambda_1}(\Omega)= \mu_{\lambda_0}(\Omega)$, then $\mu_{\lambda_0}(\Omega)$ is not achieved. 
\end{itemize}
\end{Lemma}

\proof
We take any $\lambda\ge \lambda_1$, so that $\mu_\lambda(\Omega)\le \mu_{\lambda_1}(\Omega)<\mu_{\lambda_0}(\Omega)$.
Then we 
take a minimizing sequence $(u_n)\subset H^2_0(\Omega)$ for $\mu_\lambda(\Omega)$, such that
\[
\int_\Omega|x|^{-4}|u_n|^2\,dx=1~\quad\forall n\ge 1,\quad J_{\lambda}(u_n)=\mu_{\lambda}(\Omega)+o(1)
\]
as $n\to \infty$.
Let $R>0$ be such that $B_R\supseteq \Omega$. Since
$$
\int_\Omega |\Delta u_n|^2\,dx=J_{\lambda}(u_n)+\lambda \int_\Omega |u_n|^2\,dx\le
\mu_{\lambda}(\Omega)+ |\lambda| R^{4}\int_\Omega|x|^{-4}|u_n|^2\,dx +o(1),
$$
we see that $u_n$ is bounded in $H^2_0(\Omega)$. Thanks to  Rellich's  theorem, we have that there exists $u\in H^2_0(\Omega)$ such that, up to a subsequence,
$$
\begin{gathered}
u_n\to u~~\text{weakly in $H^2_0(\Omega)$,}\quad u_n\to u~~\text{in $L^2(\Omega)$ and pointwise a.e.,}\\
|x|^{-2}u_n\to |x|^{-2}u~~\text{weakly in $L^2(\Omega)$.}
\end{gathered}
$$
Consider the sequence $\theta_n=u_n-u$. We have 
$$
\Delta\theta_n\to 0~~\text{weakly in $L^2(\Omega)$,}\quad \theta_n\to 0~~\text{in $L^2(\Omega)$,}\quad
|x|^{-2}\theta_n\to 0~~\text{weakly in $L^2(\Omega)$,}
$$
which imply
\begin{gather}
\nonumber
\int_\Omega|\Delta u_n|^2\,dx=\int_\Omega|\Delta u|^2\,dx+\int_\Omega|\Delta \theta_n|^2\,dx+o(1),\quad
\int_\Omega| u_n|^2\,dx=\int_\Omega| u|^2\,dx+o(1)\\
\nonumber
\int_\Omega|\Delta\theta_n|^2\,dx=\int_\Omega|\Delta\theta_n|^2\,dx-\lambda_0 \int_\Omega|\theta_n|^2\,dx+o(1)\ge
\mu_{\lambda_0}(\Omega) \int_\Omega|x|^{-4}|\theta_n|^2\,dx+o(1)\,,\\
 \label{thetan}
 \int_\Omega|x|^{-4}|u|^2\,dx=\int_\Omega|x|^{-4}|u_n|^2\,dx- \int_\Omega|x|^{-4}|\theta_n|^2\,dx+o(1)
 =1- \int_\Omega|x|^{-4}|\theta_n|^2\,dx+o(1).
\end{gather}
Therefore, we easily get 
$$
\begin{aligned}
\mu_{\lambda}(\Omega)&=J_{\lambda}(u_n)+o(1)=
\int_\Omega|\Delta u|^2\,dx-\lambda\int_\Omega|u|^2\,dx+\int_\Omega|\Delta \theta_n|^2\,dx+o(1)\\
&\ge\mu_{\lambda}(\Omega)\int_\Omega|x|^{-4}|u|^2\,dx+\mu_{\lambda_0}(\Omega)\int_\Omega|x|^{-4}|\theta_n|^2\,dx+o(1)
\\&
=\mu_{\lambda}(\Omega)+\left[\mu_{\lambda_0}(\Omega)-\mu_{\lambda}(\Omega)\right]\int_\Omega|x|^{-4}|\theta_n|^2\,dx+o(1).
\end{aligned}
$$
Since $\mu_{\lambda}(\Omega)<\mu_{\lambda_0}(\Omega)$, 
this implies that $|x|^{-2}\theta_n\to 0$ strongly in $L^2(\Omega)$. Finally, taking also  \eqref{thetan}
into account, we can conclude that 
$$
\int_\Omega|x|^{-4}|u|^2\,dx=1\quad\text{and}\quad J_{\lambda}(u)= \mu_{\lambda}(\Omega)\,,
$$
that is, $u$ achieves $\mu_{\lambda}(\Omega)$. This ends the proof of $i)$.

\medskip
 Assume that $\mu_{\lambda_0}(\Omega)$ is achieved by some $u\in H^2_0(\Omega)$ having unitary norm in $L^2(\Omega;|x|^{-4}dx)$.
Then 
$$
\mu_{\lambda_1}(\Omega)\le J_{\lambda_1}(u_0)=J_{\lambda_0}(u_0)-({\lambda_1}-\lambda_0)\int_\Omega|u|^2\,dx=\mu_{\lambda_0}(\Omega)-({\lambda_1}-\lambda_0)\int_\Omega|u|^2\,dx
<\mu_{\lambda_0}(\Omega)\,,
$$
hence $\mu_{\lambda_1}(\Omega)<\mu_{\lambda_0}(\Omega)$. This proves claim $ii)$.

Since  $iii)$ is an equivalent formulation of $ii)$, the lemma is completely proved.
\QED

\subsection{Proof of Theorem \ref{T:in_bordo}}

We denote by $\lambda(\Delta^{\!2},\Omega)$ the principal value of $\Delta^{\!2}$ on $H^2_0(\Omega)$, that is
$$
\lambda(\Delta^{\!2},\Omega)=\inf_{v\in H^2_0(\Omega)\atop v\neq 0}
\frac{\displaystyle\int_\Omega|\Delta v|^2\,dx}{\displaystyle\int_\Omega| v|^2\,dx}\,.
$$
By  Rellich's  theorem,  the infimum $\lambda(\Delta^{\!2},\Omega)$ is achieved. 

Evidently $\mu_{\lambda}(\Omega)=0$ if and only if $\lambda=\lambda(\Delta^2,\Omega)$. 
In particular, $\mu_{\lambda}(\Omega)>0$ for $\lambda<\lambda(\Delta^2,\Omega)$. 
Thanks to Lemma \ref{L:mu0_Omega}, this is enough to conclude that $\mu^\infty_\Omega>0$.

\medskip

To prove $i)$, take  $\lambda>\lambda^*_\Omega$. Then there exists $\lambda_0\in (\lambda^*_\Omega,\lambda)$ such that $\mu_{\lambda_0}(\Omega)$ is
achieved. Thus $\mu_\lambda(\Omega)$ is achieved as well, and $\mu_\lambda(\Omega)< \mu_{\lambda_0}(\Omega)\le \mu^\infty_\Omega$ by $ii)$ in Lemma \ref{L:compactness}.

\medskip

Take any $\lambda_0>\lambda^*_\Omega$. Since $i)$ in Theorem \ref{T:in_bordo} is proved, we can infer that $\mu_{\lambda_0}(\Omega)$ is achieved.
Therefore, \mbox{$\mu_{\lambda}(\Omega)<\mu_{\lambda_0}(\Omega)$} for any $\lambda>\lambda_0$ by $ii)$ in
Lemma \ref{L:compactness}. 
This shows that the function $\lambda\mapsto \mu_\lambda(\Omega)$ is strictly decreasing on $(\lambda^*_\Omega,\infty)$. }

Suppose that it is affine on 
$[\lambda_0,\lambda_1]\subset (\lambda^*_\Omega,\infty)$. Set $\lambda=\frac12(\lambda_0+\lambda_1)$,
so that $\frac12\left(\mu_{\lambda_0}(\Omega)+
\mu_{\lambda_1}(\Omega)\right)=\mu_{\lambda}(\Omega)$. Take 
a minimizer $u\in H^2_0(\Omega)\setminus\{0\}$ for $\mu_\lambda(\Omega)$. 
Then
\[
\mu_{\lambda}(\Omega)=J_{\lambda}(u)=\frac12\left(J_{\lambda_0}(u)+J_{\lambda_1}(u)\right)\ge \frac12\left(\mu_{\lambda_0}(\Omega)+
\mu_{\lambda_1}(\Omega)\right)=\mu_{\lambda}(\Omega).
\]
We see that $u$ achieves also $\mu_{\lambda_0}(\Omega)$ and $\mu_{\lambda_1}(\Omega)$, in addition to $\mu_\lambda(\Omega)$.
Therefore, $u$ solves
\[
\Delta^{\!2}u=\lambda_0 u+\mu_{\lambda_0}(\Omega)|x|^{-4}u\quad\text{and}\quad \Delta^{\!2}u=\lambda_1 u+\mu_{\lambda_1}(\Omega)|x|^{-4}u\quad\text{in }\Omega\,.
\]
Since $u\neq 0$, we infer that there exists a set $A\subset \Omega$ of positive measure, such that 
\[
\lambda_1- \lambda_0=( \mu_{\lambda_0}(\Omega)-\mu_{\lambda_1}(\Omega))|x|^{-4}\quad\ \text{for any $x\in A$.}
\]
This implies that $\lambda_1=\lambda_0$ and concludes the proof of $ii)$.  

\medskip

If $\lambda^*_\Omega>-\infty$, then the  function $\lambda\mapsto \mu_\lambda(\Omega)$ is constant on $(-\infty,\lambda^*_\Omega]$
because it is non increasing and thanks to $i)$ in Lemma \ref{L:compactness}. This concludes the proof.
\QED

\section{Conic domains}
\label{S:conic}

Following \cite{FM}, we say that a {\it conic domain} is the intersection of 
a cone with a ball about the origin. 
We consider here general cones, although
only the conic domains $B_R$ and $B_R\cap\mathbb R^d_+$ are needed for our purposes.

Since $\mu_\lambda(B_R\cap \mathcal C)= \mu_{R^4\lambda}(B_1\cap \mathcal C)$, we can assume, without loss of generality, 
that $R=1$. 
In the lowest dimension $d=1$, the only conic domain we will deal with is the interval $(0,1)$.

If 
$d\ge 5$ or if $\mathcal C\subsetneq \R^d\setminus\{0\}$ is a proper cone, then
the Rellich constant $\mu_0(\mathcal C)$ is positive.
Moreover, the conic domain $\Omega= B_1\cap \mathcal C$ satisfies $(A1)$ or $(A2)$, hence Theorem \ref{T:in_bordo} applies.

The main result of this section needs the notation introduced in (\ref{eq:omegaC}) and (\ref{eq:E_omega}). 

\begin{Lemma} 
\label{L:conic_domains1} 
Let $\mathcal C$ be a cone in $\R^d$; if $d\le 4$ assume that $\mathcal C$ is proper.  
Then $\mu_\lambda(B_1\cap \mathcal C)\le\mu_0(\mathcal C)$ for any $\lambda\in\R$. Moreover, 
setting
$$
\lambda^*_{B_1\cap \mathcal C}=\inf\big\{\lambda\in\R~|~\mu_\lambda(B_1\cap \mathcal C)~~\text{is achieved}\big\}\,,
$$
we have that $\mu^\infty_{B_1\cap \mathcal C}=\mu_0(\mathcal C)=\mu_\lambda(B_1\cap \mathcal C)$ for any $\lambda\le \lambda^*_{B_1\cap \mathcal C}$, and
that
the following estimates hold:
\begin{itemize}
\item[$i)$] If $d\ge 2$ and ${\rm E}(\omega_\mathcal C)+\gamma_d\ge 0$, then 
$$
\lambda^*_{B_1\cap \mathcal C}\ge 
\lambda(\Delta^{\!2},\ball)+8\left({\rm E}(\omega_\mathcal C)+\gamma_d\right)\lambda(-\Delta,\disk)\,;
$$
\item[$ii)$] If $d\ge 2$ and ${\rm E}(\omega_\mathcal C)+\gamma_d< 0$, then 
$$
\lambda^*_{B_1\cap \mathcal C}\ge \frac{{\rm E}(\omega_\mathcal C)+\gamma_d+2}{2}\lambda(\Delta^{\!2},\ball)
+\frac{|{\rm E}(\omega_\mathcal C)+\gamma_d|}{2}\lambda(-\Delta,\ball)\lambda(-\Delta,\disk)\,;
$$
\item[$iii)$] If $d=1$, then
$$\lambda^*_{(0,1)}\ge \dfrac58 \lambda(\Delta^{\!2},\ball)+\dfrac38\lambda(-\Delta,\ball)\lambda(-\Delta,\disk)\,.$$
\end{itemize}
\end{Lemma}

The proof of Lemma \ref{L:conic_domains1} needs some preliminary results.

\subsection{The Emden-Fowler transform}
\label{SS:EF}

The Emden-Fowler transform is an essential tool in \cite{CM}, on which the present paper is based.

Let $\mathcal C\subseteq \R^d$ be a  cone. If $d<4$ we assume that $\mathcal C$ is proper. If $d\ge 4$, we can take 
$\mathcal C=\R^d$, so that $B_1\cap \mathcal C = B_1$ (recall that for $d\ge 4$, $C^\infty_c(B_1\setminus\{0\})$ is dense in $H^2_0(B_1)$,
see Lemma \ref{L:density}). 
\medskip

We first consider the  case $d\ge 2$. 
We use  
spherical coordinates $r>0$, $\sigma\in \S^{d-1}$ in $\R^d$. Using \cite[Lemma 3.6]{CM}, one can show that the transform $u\mapsto \mathcal T_du$ given by
$$
u(r\sigma)=r^{\frac{4-d}{2}}(\mathcal T_du)(-\log r,\sigma)\,,
$$
is an invertible isomorphism 
$$
\mathcal T_d:H^2_0(B_1\cap \mathcal C)\to H^2_0(\R_+\times \omega_\mathcal C)\,.
$$
Setting $v=\mathcal T_du$ for $u\in H^2_0(B_1\cap \mathcal C)$, we have the equalities
\begin{equation}
\label{eq:EF}
\begin{gathered}
\int\limits_{B_1\cap \mathcal C}|x|^{-4}|u|^2\,dx=\iint\limits_{\R_+\times\omega_\mathcal C}|v|^2\,dtd\sigma~,\quad
\int\limits_{B_1\cap \mathcal C}|u|^2\,dx=\iint\limits_{\R_+\times\omega_\mathcal C} e^{-4t}|v|^2\,dtd\sigma~, \\
\int\limits_{B_1\cap \mathcal C}|\Delta u|^2\,dx=\iint\limits_{\R_+\times\omega_\mathcal C}|\Delta_\sigma v-\gamma_dv|^2\,dtd\sigma+
\underbrace{\iint\limits_{\R_+\times\omega_\mathcal C}\big[|v_{tt}|^2+2|\nabla_\sigma v_t|^2+2(\gamma_d+2)|v_t|^2\big]\,dtd\sigma}_{\displaystyle =:\,G(v)}\,,
\end{gathered}
\end{equation}
so that 
\begin{equation}
\label{eq:EF_2a}
\mu_\lambda(B_1\cap \mathcal C)=\inf_{v\in H^2_0(\R_+\times \omega_\mathcal C)\atop
v\neq 0}
\left[~\frac{\displaystyle\iint\limits_{\R_+\times\omega_\mathcal C}|\Delta_\sigma v-\gamma_dv|^2\,dtd\sigma}{\displaystyle\iint\limits_{\R_+\times\omega_\mathcal C}|v|^2\,dtd\sigma}
+\frac{G(v)-\lambda \displaystyle\iint\limits_{\R_+\times\omega_\mathcal C}e^{-4t}|v|^2\,dtd\sigma}{\displaystyle\iint\limits_{\R_+\times\omega_\mathcal C}|v|^2\,dtd\sigma}~\right].
\end{equation}
Moreover, the minimization problems in (\ref{eq:EF_2a})  and in (\ref{eq:mu_def}), with $\Omega=B_1\cap\mathcal C$, are equivalent. 

In the lowest dimensional case $d=1$ we let $\mathcal C=\R_+$; the Emden-Fowler transform is defined by
$$
\mathcal T_1:H^2_0(0,1)\to H^2_0(\R_+)~,\quad u(x)=r^{\frac{3}{2}}(\mathcal T_du)(-\log x)\,.
$$
Formula (\ref{eq:EF_2a}) then reduces to the corresponding one-dimensional expression upon suppressing integration and derivatives with respect to the 
$\sigma$-variable. Since $\gamma_1=-3/4$, we now obtain that
\begin{equation}
\label{eq:EF_Unoa}
\mu_\lambda(0,1)=\frac{9}{16}+\inf_{v\in H^2_0(\R_+)\atop
v\neq 0} 
\frac{\displaystyle\int_0^\infty \big[|v_{tt}|^2+\tfrac52|v_t|^2\big]\,dt
-\lambda \displaystyle \int_0^\infty e^{-4t}|v|^2\,dt}{\displaystyle\int_0^\infty|v|^2\,dt}\,.
\end{equation}

\subsection{Auxiliary best constants}
\label{SS:Auxiliary}

In order to face the cases $d\ge 2$ and $d=1$ in a unified way, we study the best constant
\begin{equation}
\label{eq:Gamma_def}
\Lambda_\#(\beta) =\inf_{f\in H^2_0(\R_+)\atop f\neq 0} \frac{\displaystyle{\int_0^\infty\big[|f_{tt}|^2+2\beta|f_t|^2\big]dt}}{\displaystyle\int_0^\infty e^{-4t}|f|^2 dt}\,,
\quad \beta\ge 0\,,
\end{equation}
which is related to the infima in (\ref{eq:EF_2a}), (\ref{eq:EF_Unoa}) for $\beta=\gamma_d+2$.
Using two times the Hardy inequality (\ref{eq:Hardy1}),  one  easily obtains
$$
\Lambda_\#(\beta) \ge \inf_{f\in H^2_0(\R_+)\atop f\neq 0} \frac{\displaystyle{\frac{9}{16}\int_0^\infty t^{-4}|f|^2\,dt+2\beta\frac{1}{4}\int_0^\infty t^{-2}|f|^2\,dt}}{\displaystyle\int_0^\infty e^{-4t}|f|^2 dt}\ge \frac{9}{16}e^4+2\beta e^2.
$$
The proof of Lemma \ref{L:conic_domains1} needs a more refined estimate. 

\begin{Lemma}
\label{L:diesis_achieved}
Let $\beta\ge 0$. Then
\begin{equation}
\label{eq:lambda_beta_positive}
\Lambda_\#(\beta)\ge 
\begin{cases}
~\lambda(\Delta^{\!2},\ball)+8(\beta-2)\lambda(-\Delta,\disk)&\text{if $\beta\ge 2$,}\\ \\
~\dfrac\beta2 \lambda(\Delta^{\!2},\ball)+\dfrac{2-\beta}{2}\lambda(-\Delta,\ball)\lambda(-\Delta,\disk)&\text{if $0\le \beta< 2$.}
\end{cases}
\end{equation}
\end{Lemma}

\proof
We  reformulate (\ref{eq:Gamma_def}) as a minimization problem 
over radial functions in $H^2_0(\ball)$, via the Emden--Fowler transform
$$
(\mathcal T_4w)(x)= f(-\log |x|)\,,\qquad w\in H^2_{0, \rm rad}(\ball):=\big\{ w\in H^2_0(\ball)~|~ w=w(|x|)~\big\}\,.
$$
It turns out that  $\mathcal T_4$ is  an isomorphism $H^2_{0, \rm rad}(\ball)\to  H^2_0(\R_+)$, in view of Lemma \ref{L:density}.

Let $w\in H^2_{0, \rm rad}(\ball)$.
Since $f=\mathcal T_4 w$ does not depend on $\sigma\in\S^3$ and since $\gamma_4=0$, by applying (\ref{eq:EF})
with $\mathcal C=\ball$ and $\omega_\mathcal C=\S^3$, we get
$$
\begin{gathered}
\int_{\ball}|\Delta w|^2~\!dx=|\S^3|\int_0^\infty\big(|f_{tt}|^2+4|f_t|^2\big)\,dt~,\quad
\int_{\ball}|w|^2~\!dx=|\S^3|\int_0^\infty e^{-4t}|f|^2\,dt\,.
\end{gathered}
$$
Since in addition 
$$
\int_{\ball}|x|^{-2}|\nabla w|^2~\!dx=|\S^3|\int_0^\infty|f_t|^2\,dt\,,
$$
we infer that 
$$
\Lambda_\#(\beta)=\inf_{w\in H^2_{0,\rm rad}(\ball)\atop w\neq 0} 
\frac{\displaystyle{\int_{\ball}|\Delta w|^2~\!dx+2(\beta-2)\int_{\ball}|x|^{-2}|\nabla w|^2~\!dx}}{\displaystyle\int_{\ball}|w|^2 dx}\,.
$$
If $\beta\ge 2$, then inequality (\ref{eq:ine_1})  in  Appendix \ref{A:BV} implies
$$
\Lambda_\#(\beta)\ge\inf_{w\in H^2_{0,\rm rad}(\ball)\atop w\neq 0} 
\frac{\displaystyle{\int_{\ball}|\Delta w|^2~\!dx+8(\beta-2)\lambda(-\Delta,\disk)\int_{\ball}|w|^2~\!dx}}{\displaystyle\int_{\ball}|w|^2 dx}
=\lambda(\Delta^{\!2},\ball)+8(\beta-2)\lambda(-\Delta,\disk)\,.
$$
If $\beta\in[0,2)$, we use the inequality (\ref{eq:ine_2}) to estimate
$$
\begin{aligned}
\Lambda_\#(\beta)&\ge\inf_{w\in H^2_{0,\rm rad}(\ball)\atop w\neq 0} 
\frac{\displaystyle{\frac\beta2\int_{\ball}|\Delta w|^2~\!dx+\frac{2-\beta}{2}\lambda(-\Delta,\disk)\int_{\ball}|\nabla w|^2~\!dx}}{\displaystyle\int_{\ball}|w|^2 dx}
\\ &\ge~ \frac\beta2\lambda(\Delta^{\!2},\ball)+\frac{2-\beta}{2}\lambda(-\Delta,\disk)\lambda(-\Delta,\ball).
\end{aligned}
$$
The lemma is completely proved.
\QED

\begin{Remark}
\label{R:lambda_diesis}
For completeness, we observe that the infimum $\Lambda_\#(\beta)$ is achieved because 
$H^2_0(\R_+)$ is compactly embedded into $L^2(\R_+;e^{-4t}dt)$. 
This is a consequence of  Rellich's  theorem, which gives the compactness of the embedding
$H^2_0(\ball)\hookrightarrow L^2(\ball)$, via  the isomorphisms $\mathcal T_d:H^2_{0,\rm rad}(\ball)\to H^2_0(\R_+)$, 
and $\mathcal T_d^{-1}: L^2(\R_+;e^{-4t}dt)\to L^2(\ball)$.
\end{Remark}

\medskip

\begin{Lemma}
\label{L:diesis}
Let $\omega$ be a domain in $\S^{d-1}$, $d\ge 2$. Then  
$$
M_d(\omega):=\inf_{v\in H^2_0(\R_+\times \omega)\atop v\neq 0} 
\frac{\displaystyle{\int_0^\infty \hskip-0.2cm dt \int_\omega\big[|v_{tt}|^2+2|\nabla_{\!\sigma}v_t|^2+2({\gamma}_{d}+2)|v_t|^2\big]d\sigma}}
{\displaystyle{\int_0^\infty \hskip-0.2cm e^{-4t} dt\int_\omega|v|^2 d\sigma}}\ge \Lambda_\#({\rm E}(\omega)+{\gamma}_{d}+2)\,.
$$
\end{Lemma}

\proof
Let $v\in H^2_0(\R_+\times\omega)$.
Recall that $v_t(\cdot,t)\in H^1_0(\omega)$ for a.e. $t\in\R_+$, and that
$v_t(x,\cdot)\in H^1_0(\R_+)$ for a.e. $x\in\omega$. Therefore, by the
definition of ${\rm E}(\omega)$ in (\ref{eq:m_omega}) and by (\ref{eq:Gamma_def}),
 we can estimate
$$
\begin{aligned}
\int_0^\infty \hskip-0.2cm dt \int_\omega\big[|v_{tt}|^2+2|\nabla_{\!\sigma}v_t|^2+2({\gamma}_{d}+2)&|v_t|^2\big]d\sigma
\ge \int_\omega d\sigma\int_0^\infty  \big[|v_{tt}|^2+2({\rm E}(\omega)+{\gamma}_{d}+2)|v_t|^2\big]dt\\
&\ge \Lambda_\#({\rm E}(\omega)+{\gamma}_{d}+2) \int_0^\infty \hskip-0.2cm e^{-4t} dt\int_\omega|v|^2 d\sigma\,,
\end{aligned}
$$
which concludes the proof.
\QED

\begin{Remark}
For completeness, we observe that $M_d(\omega)= \Lambda_\#({\rm E}(\omega)+{\gamma}_{d}+2)$. In fact, 
by Remark \ref{R:lambda_diesis} we can take $f\in H^2_0(\R_+)$ achieving
$\Lambda_\#({\rm E}(\omega)+\gamma_d+2)$. Then we take a 
sequence $(\f_n)_n$ of functions in $C^\infty_c(\omega)$ such that $\f_n\to \f_1$ in
$H^1_0({\omega_{\mathcal C}})$, where $\f_1\in H^2({\omega_{\mathcal C}})\cap H^1_0({\omega_{\mathcal C}})$ is the first 
eigenfunction for $-\Delta_\sigma$ in $H^1_0({\omega_{\mathcal C}})$. 
One can test $M_d(\omega)$ with $f\f_n\in H^2_0(\R_+\times \omega)$ to get, in the limit as $n\to \infty$, the opposite inequality 
$M_d(\omega)\le \Lambda_\#({\rm E}(\omega)+{\gamma}_{d}+2)$. 

Finally,  notice that from the equality $M_d(\omega)= \Lambda_\#({\rm E}(\omega)+{\gamma}_{d}+2)$ it follows that the infimum $M_d(\omega)$ is not achieved.
\end{Remark}

\subsection{Proof of Lemma \ref{L:conic_domains1}}
We fix any nontrivial function $u\in C^\infty_c(\mathcal C)$ and put 
$u_h(x)=u(hx)$ for $h>0$. If $h$ is large enough, then $u_h\in C^\infty_c(B_1\cap \mathcal C)$, and therefore
for any $\lambda\in \R$ we can estimate
$$
\mu_\lambda(B_1\cap \mathcal C)\le \frac
{\displaystyle\int_\mathcal C|\Delta u_h|^2~\!dx-\lambda\displaystyle\int_\mathcal C|u_h|^2~\!dx}
{\displaystyle\int_\mathcal C|x|^{-4}|u_h|^2~\!dx}
=\frac
{h^4\displaystyle\int_\mathcal C|\Delta u|^2~\!dx-\lambda\displaystyle\int_\mathcal C|u|^2~\!dx}
{h^4\displaystyle\int_\mathcal C|x|^{-4}|u|^2~\!dx}=
\frac
{\displaystyle\int_\mathcal C|\Delta u|^2~\!dx}
{\displaystyle\int_\mathcal C|x|^{-4}|u|^2~\!dx}+o(1)\,,
$$
as $h\to \infty$. Since $u$ was arbitrarily chosen, this implies that 
$$
\mu_\lambda(B_1\cap \mathcal C)\le \inf_{u\in C^\infty_c(\mathcal C)\atop u\neq 0} 
\frac
{\displaystyle\int_\mathcal C|\Delta u|^2~\!dx}
{\displaystyle\int_\mathcal C|x|^{-4}|u|^2~\!dx} =\mu_0(\mathcal C)\,.
$$

Easily, for any $\lambda\le 0$ we have $\mu_\lambda(B_1\cap \mathcal C)\ge \mu_0(B_1\cap\mathcal C)\ge \mu_0(\mathcal C)$.
We proved that $\mu_\lambda(B_1\cap \mathcal C)=  \mu_0(\mathcal C)$ for any $\lambda\le 0$.
Therefore, $\lambda^*_{B_1\cap \mathcal C}$ is nonnegative, and
$\mu^\infty_\Omega=\mu_0(\mathcal C)=\mu_\lambda(B_1\cap \mathcal C)$ for any $\lambda\le \lambda^*_{B_1\cap \mathcal C}$
by Theorem \ref{T:in_bordo}.

\medskip

We now assume that $d\ge 2$ and prove the estimates in $i)-ii)$. 
Take any $\lambda> \lambda^*_{B_1\cap \mathcal C}$. Then
\begin{equation}
\label{eq:strict}
\mu_{\lambda}(B_1\cap \mathcal C)<\mu_0(\mathcal C)
\end{equation}
and $\mu_{\lambda}(B_1\cap \mathcal C)$ is attained by Theorem \ref{T:in_bordo}.
 Take a minimizer  $u^\lambda\in H^2_0(B_1\cap \mathcal C)\setminus\{0\}$ for $\mu_{\lambda}(B_1\cap \mathcal C)$,
and put 
$$v^\lambda:=\mathcal T_d u^\lambda\,,
$$
where
$\mathcal T_d$ is the Emden-Fowler transform in Subsection \ref{SS:EF}. 
We have that
$v^\lambda\in H^2_0(\R_+\times\omega_\mathcal C)$ achieves the infimum in (\ref{eq:EF_2a}). 
We can assume that $v^\lambda$ has unitary norm in $L^2(\R_+\times\omega_\mathcal C)$, so that 
$$
\begin{aligned}
\mu_{\lambda}(B_1\cap\mathcal C)&=
\iint\limits_{\R_+\times\omega_\mathcal C}|\Delta_\sigma v^\lambda-\gamma_dv^\lambda|^2\,dtd\sigma
\\&+
\iint\limits_{\R_+\times \omega_\mathcal C}\big[|v^\lambda_{tt}|^2+2|\nabla_\sigma v^\lambda_t|^2+2(\gamma_d+2)|v^\lambda_t|^2\big]\,dtd\sigma-
\lambda \iint\limits_{\R_+\times\omega_\mathcal C}e^{-4t}|v^\lambda|^2\,dtd\sigma\,.
\end{aligned}
$$
By the definition of $m_{\gamma_d}(\omega_\mathcal C)$ in (\ref{eq:m_omega}), thanks to (\ref{eq:mu_m}), and using Lemma \ref{L:diesis}, we can estimate
$$
\begin{aligned}
\iint\limits_{\R_+\times\omega_\mathcal C}|\Delta_\sigma v^\lambda-\gamma_dv^\lambda|^2\,dtd\sigma&\ge m_{\gamma_d}(\omega_\mathcal C)=\mu_0(\mathcal C)\,,
\\
\iint\limits_{\R_+\times \omega_\mathcal C}\big[|v^\lambda_{tt}|^2+2|\nabla_\sigma v^\lambda_t|^2+2(\gamma_d+2)|v^\lambda_t|^2\big]\,dtd\sigma
&\ge \Lambda_\#({\rm E}(\omega_\mathcal C)+\gamma_d+2) \iint\limits_{\R_+\times\omega_\mathcal C}e^{-4t}|v^\lambda|^2\,dtd\sigma\,.
\end{aligned}
$$
Taking also (\ref{eq:strict}) into account, we  infer that
$$
0>\mu_{\lambda}(B_1\cap\mathcal C)- \mu_0(\mathcal C)
\ge 
\big[\Lambda_\#({\rm E}(\omega_\mathcal C)+\gamma_d+2)-\lambda\big]
\iint\limits_{\R_+\times\omega_\mathcal C}e^{-4t}|v^\lambda|^2\,dtd\sigma\,.
$$
Thus $\lambda>\Lambda_\#({\rm E}(\omega_\mathcal C)+\gamma_d+2\big)$. Since $\lambda>\lambda^*_{B_1\cap \mathcal C}$
was arbitrarily chosen, we proved that $\lambda^*_{B_1\cap \mathcal C} \ge \Lambda_\#({\rm E}(\omega_\mathcal C)+\gamma_d+2\big)$. 
Then $i)$ and $ii)$ follow from Lemma \ref{L:diesis_achieved}, with $\beta= {\rm E}(\omega)+{\gamma}_{d}+2$. 

\medskip

If $d=1$ and the conic domain is $(0,1)$, the above argument can be easily adapted. In this case the conclusion follows
by using (\ref{eq:EF_Unoa}) and the second inequality in  (\ref{eq:lambda_beta_positive}), with $\beta=5/4$.
\QED

\paragraph{Proof of Theorem \ref{T:in_Omega}}
Since $\mu_{\lambda}(\Omega)= \mu_{R^{4}\lambda}(R^{-1}\Omega)$, we can assume that $\Omega\subseteq B_1$.

Let $\delta>0$ be such that $B_\delta\subset \Omega$. Then 
$\mu_\lambda(B_1)\le \mu_\lambda(\Omega)\le \mu_\lambda(B_\delta)=\mu_{\delta^{-4}\lambda}(B_1)$
for any $\lambda\in\R$. In particular, $\mu_\lambda(\Omega)=\mu_0(\R^d)$ for any $\lambda\le \lambda^*_{B_1}=\delta^{-4}\lambda^*_{B_\delta}$ by Lemma 
\ref{L:conic_domains1}. Thanks to Theorem \ref{T:in_bordo}, we infer that $\mu^\infty_\Omega=\mu_0(\R^d)$
and that $\lambda^*_\Omega\ge \lambda^*_{B_1}$.

To conclude the proof, notice that $\lambda^*_{B_1}\ge \lambda(\Delta^{\!2},\ball)+2(d-4)\lambda(-\Delta,\disk)$
by $i)$ in Lemma \ref{L:conic_domains1} (with 
$\mathcal C=\R^d$ and $\gamma_d=d(d-4)/4>0$). 
\QED

\paragraph{Proof of Theorem \ref{T:Uno}} 
Since $\mu_{\lambda}(0,R)= \mu_{R^{4}\lambda}(0,1)$, we can assume that $R=1$. 
Using also a rescaling argument it is easy to show that $\mu_0(0,1)=\mu_0(\R_+)=9/16$, compare with (\ref{eq:916}).
Hence $\mu^\infty_{(0,R)}=9/16$ by $iv)$ in Theorem \ref{T:in_bordo}.
Finally, the estimate (\ref{eq:estimate*_Uno}) follows from $iii)$ in Lemma \ref{L:conic_domains1}.
\QED

\section{Proof of Theorem \ref{T:smooth}}
\label{S:proofs}

For any $\eps>0$ there exists 
a concave cone $\mathcal C^+_\eps$, a convex cone $\mathcal C^-_\eps$, and $\delta>0$ such that
\begin{equation}
\label{eq:eps_cones}
\begin{gathered}
0<\mu_0(\R^d_+)-\mu_0(\mathcal C^+_\eps)<\eps\quad \text{and} \quad \Omega\cap B_\delta\subset \mathcal C^+_\eps\,,\\
0<\mu_0(\mathcal C^-_\eps)-\mu_0(\R^d_+)<\eps\quad \text{and} \quad B_\delta\cap \mathcal C^-_\eps\subset \Omega\,.
\end{gathered}
\end{equation}
This elementary fact can be proved as follows. Take 
$\delta_0>0$ and a continuous 
function $F\colon {B}_{\delta_0}\to\R$, differentiable at $0$ 
and satisfying (\ref{eq:F_bordo}).
By continuity, we have that 
\[
\lim_{\rho\searrow 0}\mu_0(C_\rho)=\mu_0(\R^d_+)~,\text{where}\quad 
C_\rho=\{x\in\R^d~\colon~ x\cdot\nabla F(0)<\rho\left|x\right|\,\}\supset \R^d_+\,.
\]
Hence $\mu_0(C_{\rho_\eps})>\mu_0(\R^d_+)-\eps$ if $\rho_\eps>0$ is small. Set $\mathcal C^+_\eps:=C_{\rho_\eps}$.
Since $0>F(x)=x\cdot \nabla F(0)+o(|x|)$ as $x\in\Omega, x\to 0$, we have that
if $\delta\le \delta_0$ is small enough and $x\in \Omega\cap B_\delta$, then $x\in\mathcal C^+_\eps$.  
This proves the existence of the concave  cone $\mathcal C^+_\eps$. For $\mathcal C^-_\eps$
argue in a similar way.  

\medskip

We are in position to prove (\ref{eq:limit}). In fact, it is equivalent to establish that
\begin{equation}
\label{eq:limit2}
\lim\limits_{\lambda\to -\infty}\mu_\lambda(\Omega)=\mu_0(\R^d_+)\,,
\end{equation}
use $iii)$ in Theorem \ref{T:in_bordo} in case $\lambda^*_\Omega$ is finite. 

Take any $\eps>0$ and find the cones $\mathcal C^\pm_\eps$ and $\delta>0$ as in (\ref{eq:eps_cones}). 
Thanks to Lemma \ref{L:conic_domains1} and by
the monotonicity of $\mu_\lambda$ with respect to set inclusion, 
for any $\lambda\le 0$ we have
$\mu_0(\mathcal C^+_\eps)=\mu_0(B_\delta\cap\mathcal C^+_\eps)\le \mu_0(\Omega\cap B_\delta)$ and
$\mu_\lambda(\Omega)\le \mu_\lambda(B_\delta\cap \mathcal C^-_\eps)=\mu_0(\mathcal C^-_\eps)$ for any 
$\lambda\le 0$.
We infer that 
\begin{gather}
\label{eq:step1a}
\mu_0(\R^d_+)-\eps<\mu_0(\Omega\cap B_\delta)\\
\label{eq:step1b}
\mu_\lambda(\Omega)<\mu_0(\R^d_+)+\eps\quad \text{for any $\lambda\le 0$} \,.
\end{gather}
Now we take a function $\eta\in C^\infty(\R)$ such that 
$0\le \eta\le 1$, $\eta\equiv 0$ on $(-\infty,1)$ and $\eta\equiv 1$ on $(2,\infty)$. Let
\begin{equation}
\label{eq:psih}
\psi_h(x)=\eta\left(\frac{-\log\left|x\right|}h\right)=
\begin{cases}
0&\text{ if $\left|x\right|\ge e^{-h}$}\\ 
1&\text{ if $\left|x\right|\le e^{-2h}$.}
\end{cases}
\end{equation}
Then $\psi_h\in C^\infty_c(\R^d)$, $\text{supp}\nabla\psi_h\subseteq \{e^{-2h}\le|x|\le e^{-h}\}$, and 
there exists a constant $c>0$, depending only on $\eta$, such that
\begin{equation}
\label{eq:psi_derivative}
|\nabla \psi_h(x)|\le c h^{-1}|x|^{-1}~,\quad |\Delta\psi_h(x)|\le c h^{-1}|x|^{-2}\qquad \text{if~~ $e^{-2h}<|x|<e^{-h}$}\,.
\end{equation}
Further, for $h>|\log\delta|$  and for any $u\in H^2_0(\Omega)$ we have that $\psi_h u\in H^2_0(\Omega\cap B_\delta)$, and therefore 
 (\ref{eq:step1a}) gives
\begin{equation}
\label{eq:piove}
\big(\mu_0(\R^d_+)-\eps\big)\int_\Omega|x|^{-4}|\psi_h u|^2\,dx 
\le \int_\Omega|\Delta(\psi_hu)|^2\,dx.
\end{equation}
We write $u=\psi_hu+(1-\psi_h)u$ to estimate
$$
\begin{aligned}
\int_{\Omega}|x|^{-4}|u|^2\,dx&=\int_{\Omega}|x|^{-4}|\psi_hu|^2\,dx+2 \int_{\Omega}|x|^{-4}\psi_h(1-\psi_h)|u|^2\,dx
+\int_{\Omega}|x|^{-4}(1-\psi_h)^2|u|^2\,dx\\
&\le \int_{\Omega}|x|^{-4}|\psi_hu|^2\,dx+3e^{8h}\int_\Omega|u|^2\,dx.
\end{aligned}
$$
Thanks to (\ref{eq:psi_derivative}) and since  $\mu_0^{\rm HR}(\Omega)>0$, $\mu_0(\Omega)>0$
by Lemmata \ref{L:mu0_Omega0} and \ref{L:mu0_Omega}, we obtain
$$
\begin{gathered}
\int_\Omega|\nabla \psi_h\cdot\nabla u|^2\,dx\le {c}{h^{-2}}\int_\Omega|x|^{-2}|\nabla u|^2\,dx
\le {c}{h^{-2}} \int_\Omega|\Delta u|^2\,dx\,,\\
\int_\Omega|u\Delta  \psi_h|^2\,dx \le {c}{h^{-2}}\int_\Omega|x|^{-4}| u|^2\,dx
\le {c}{h^{-2}} \int_\Omega|\Delta u|^2\,dx\,,
\end{gathered}
$$
where the constants $c$ depend on $\eta$, $\Omega$ but not on $u$.
Therefore, the triangle inequality gives 
$$
\int_\Omega|\Delta(\psi_h u)|^2\,dx\le\left( \|\Delta u\|_2+2\|\nabla \psi_h\cdot\nabla u\|_2+\|u\Delta \psi_h\|_2\right)^2
\le  \left(1+ch^{-1}\right)^2\int_\Omega|\Delta u|^2\,dx.
$$
Collecting the above estimates, from (\ref{eq:piove}) we infer that 
$$
\big(\mu_0(\R^d_+)-\eps\big)\Big(\int_{\Omega}|x|^{-4}|u|^2\,dx-3e^{8h}\int_\Omega|u|^2\,dx\Big)
\le  \left(1+ch^{-1}\right)^2\int_\Omega|\Delta u|^2\,dx\,,
$$
which gives
$$
\frac{\mu_0(\R^d_+)-\eps}{\left(1+ch^{-1}\right)^2}\int_{\Omega}|x|^{-4}|u|^2\,dx\le \int_\Omega|\Delta u|^2\,dx-\lambda_h\int_\Omega|u|^2\,dx~,
\quad \text{where} ~  \lambda_h=-\frac{3\big(\mu_0(\R^d_+)-\eps\big)e^{8h}}{(1+ch^{-1})^2}.
$$
Notice that $\lambda_h\to -\infty$ as $h\to \infty$. Since $u\in H^2_0(\Omega)$ was arbitrarity chosen, we have proved that
$$
\mu_{\lambda_h}(\Omega)\ge \frac{\mu_0(\R^d_+)-\eps}{\left(1+ch^{-1}\right)^2}= \mu_0(\R^d_+)-\eps+o(1)\quad
\text{as $h\to\infty$.}
$$
Therefore, taking also the monotonicity of $\lambda\mapsto \mu_\lambda(\Omega)$ and (\ref{eq:step1b}) into account, we have that
$$
\mu_0(\R^d_+)-\eps \le \lim_{h\to \infty}\mu_{\lambda_h}(\Omega)=\lim_{\lambda\to-\infty}\mu_\lambda(\Omega)\le \mu_0(\R^d_+)+\eps.
$$
Since $\eps>0$ was arbitrarily chosen, then (\ref{eq:limit2}) follows.

\medskip

It only remains to prove the estimate in  (\ref{eq:estimate_half}) in case $\Omega\subset\R^d_+$. 
Likely in the proof of Theorem \ref{T:in_Omega}, 
we can assume that $\Omega\subset B_1\cap \R^d_+$. Further, we have $\mu^\infty_\Omega=\mu^\infty_{B_1\cap \R^d_+}=\mu_0(\R^d_+)$
by (\ref{eq:limit}) and Lemma \ref{L:conic_domains1}, and $\mu_\lambda(B_1\cap \R^d_+)\le \mu_\lambda(\Omega)$ for any $\lambda\in\R$.
Therefore, from Theorem \ref{T:in_bordo} we infer that $\lambda^*_{B_1\cap \R^d_+}\le \lambda^*_\Omega$, 
and we are led to estimate $\lambda^*_{B_1\cap \R^d_+}$ from below.

Now, we take $\mathcal C=\R^d_+$ in Lemma \ref{L:conic_domains1}, so that $\omega_\mathcal C=\S^{d-1}_+$
is a hemisphere. Since 
$$
{\rm E}(\S^{d-1}_+)+\gamma_d=(d-1)+\frac{d(d-4)}{4}=\frac{d^2-4}{4}\ge 0\,,
$$
then (\ref{eq:estimate_half}) for $\Omega=B_1\cap \R^d_+$ follows from Lemma \ref{L:conic_domains1}. The proof is complete.
\QED

\appendix

\section{Examples} 
\label{A:Example} 

\begin{Example}
\label{Example}
Let $\delta\in (0,1)$, $\nu\in\S^{d-1}$, and consider the concave cone $\mathcal C_\delta=\{ x\in\R^N~|~x\cdot\nu>-\delta|x|~\}$.
There exists $R_\delta>1$ such that for any domain $\Omega\subset\R^d$ satisfying the assumptions in Theorem \ref{T:smooth} and such that
$$
\Omega\supseteq \mathcal C_\delta\cap \{1<|x|<R_\delta\}\,,
$$
it holds that $\lambda^*_\Omega<0$. In particular, the Rellich constant $\mu_0(\Omega)$ is achieved.
\end{Example}

This is an adaptation of the example in \cite[Proposition 6.3]{FM}. Notice that the curvature of $\Omega$ at $0\in\partial\Omega$
does not affect the conclusion.

Since
the cone $\mathcal C_\delta$ contains a half space, then its Rellich constant $\mu_0(\mathcal C_\delta)$ is smaller than $\mu_0(\R^d_+)$. Thus
there exists $u_\delta\in C^\infty_c(\mathcal C_\delta)$ such that
$$
\int_{\mathcal C_\delta}|x|^{-2}|u_\delta|^2~dx=1~,\qquad \int_{\mathcal C_\delta}|\Delta u_\delta|^2~dx  <\mu_0(\R^d_+) \,.
$$
Up to rescaling, we can assume that the support of $u_\delta$ is contained in an annulus of internal radius $1$ and
external radius $R_\delta>1$. 

If the domain $\Omega$ contains the set $\mathcal C_\delta\cap \{1<|x|<R_\delta\}$, 
then $u_\delta$ has compact support in $\Omega$. Therefore,  testing $\mu_0(\Omega)$ with $u_\delta$
we plainly get
$$
\mu_0(\Omega)\le \int_{\mathcal C_\delta}|\Delta u_\delta|^2~dx< \mu_0(\R^d_+)=\mu^\infty_\Omega\,.
$$

\begin{Example}
\label{Example2}
Set
$$
\Omega=\big\{x=(x',x_d)\in\R^{d-1}\times\R~|~ 0<x_d<|x'|^2\,\big\}\cap B_1\subset\R^d_+\,.
$$
Then $\mu^\infty_\Omega=\infty$. Therefore, $\mu_\lambda(\Omega)$
is achieved for any $\lambda\in\R$ by Theorem \ref{T:in_bordo}, and $\lambda^*_\Omega=-\infty$. 

In particular, the Rellich constant
$\mu_0(\Omega)$ is achieved. 
\end{Example}

\medskip

We start by noticing that $\Omega$ satisfies $(A2)$, hence Theorem \ref{T:in_bordo} applies. 
To prove that $\mu^\infty_\Omega=\infty$, for any small  $\theta>0$,  we set
$$
\mathcal C_\theta 
=\big\{x=(x',x_d)\in\R^{d-1}\times\R~|~ 0<x_d <\theta|x|\,\big\}~,\quad
\omega_\theta =\mathcal C_\theta\cap \S^{d-1}=\big\{\sigma\in \S^{d-1}~|~ 0<\sigma\cdot\rme<\theta\,\,\big\}\,,
$$
where $\rme=(0,\dots,0,1)$. Notice that 
$$
\mu_0(\Omega\cap B_{\theta})\ge \mu_0(\mathcal C_\theta)\,,
$$
because $\Omega\cap B_\theta\subset \mathcal C_\theta$.
We claim that
\begin{equation}
\label{eq:limCtheta}
\lim_{\theta\to 0}\mu_0(\mathcal C_\theta)=\infty\,.
\end{equation}
We prove (\ref{eq:limCtheta}) by exploiting the equality (\ref{eq:mu_m}), which gives
$$
\mu_0(\mathcal C_\theta)=
\inf_{\f\in C^\infty_c(\omega_\theta)\atop \f\neq 0}
\frac{\displaystyle{\int_{\omega_\theta}|\Delta_\sigma\f-\gamma_d\f|^2\,d\sigma}}{\displaystyle\int_{\omega_\theta}|\f|^2\,d\sigma}
=\gamma_d^2+
\inf_{\f\in C^\infty_c(\omega_\theta)\atop \f\neq 0}
\frac{\displaystyle{\int_{\omega_\theta}|\Delta_\sigma\f|^2+2\gamma_d\int_{\omega_\theta}|\nabla_{\!\sigma}\f|^2\,d\sigma}}{\displaystyle\int_{\omega_\theta}|\f|^2\,d\sigma}\,.
$$
Let ${\rm E}(\omega_\theta)$ be the first eigenvalue of $-\Delta_\sigma$ on $H^1_0(\omega_\theta)$, compare with (\ref{eq:E_omega}). 
It holds that  
$$
\lim_{\theta\to 0}{\rm E}(\omega_\theta)=\infty\,.
$$
(Indeed, ${\rm E}(\omega_\theta)=O(\theta^{-2})$, see for instance \cite[Theorem 3.3]{Bor}). 
Since for any  $\f\in C^\infty_c(\omega_\theta)$ we have
$$
\begin{aligned}
\int_{\omega_\theta}|\nabla_{\!\sigma}\f|^2\,d\sigma&=-\int_{\omega_\theta}\f\Delta_\sigma\f\,d\sigma
\le \frac{1}{\sqrt{{\rm E}(\omega_\theta)}}\Big(\int_{\omega_\theta}|\nabla_{\!\sigma}\f|^2\,d\sigma\Big)^\frac12
\Big(\int_{\omega_\theta}|\Delta_\sigma\f|^2\,d\sigma\Big)^\frac12\,,
\end{aligned}
$$
it follows that for $\theta>0$ small enough (which is needed only if $d=2,3$), we can estimate
$$
\begin{aligned}
\mu_0(\mathcal C_\theta)&\ge 
\gamma_d^2+
\big({\rm E}(\omega_\theta)+2\gamma_d\big)\inf_{\f\in C^\infty_c(\omega_\theta)\atop \f\neq 0}
\frac{\displaystyle{\int_{\omega_\theta}|\nabla_{\!\sigma}\f|^2\,d\sigma}}{\displaystyle\int_{\omega_\theta}|\f|^2\,d\sigma}
\ge  
\big({\rm E}(\omega_\theta)+2\gamma_d\big){\rm E}(\omega_\theta)\,,
\end{aligned}
$$
and (\ref{eq:limCtheta}) follows.

\medskip

To conclude the proof of the equality $\mu^\infty_\Omega=\infty$, one only needs to adapt some of the arguments in the proof of (\ref{eq:limit2}). 
We add few details for the convenience of the reader.

For any integer $h>|\log \theta|$, let $\psi_h\in C^\infty_c(B_{\theta})$ be the cut-function in (\ref{eq:psih}). 
Take any nontrivial function $u\in H^2_0(\Omega)$. 
We can  test $\mu_0 (\Omega\cap B_\theta)$ with $\psi_h u\in H^2_0(\Omega\cap B_\theta)$. Since $\mu_0 (\Omega\cap B_\theta)\ge\mu_0(\mathcal C_\theta)$, we obtain 
$$
\mu_0(\mathcal C_\theta)\int_\Omega|x|^{-4}|\psi_h u|^2\,dx\le \int_\Omega|\Delta(\psi_hu)|^2\,dx.
$$
The estimates in the proof of Theorem \ref{T:smooth} give
$$
\mu_0(\mathcal C_\theta)\Big(\int_{\Omega}|x|^{-4}|u|^2\,dx-3e^{8h}\int_\Omega|u|^2\,dx\Big)
\le \left(1+ch^{-1}\right)^2\int_\Omega|\Delta u|^2\,dx\,.
$$
Since $u$ was arbitrarily chosen in $H^2_0(\Omega)$, one infers that 
$$
\mu_{\lambda_h}(\Omega)\ge \frac{\mu_0(\mathcal C_\theta)}{\left(1+ch^{-1}\right)^2}
~,\quad\text{where}\quad 
\lambda_h= -\frac{3\mu_0(\mathcal C_\theta)e^{8h}}{(1+ch^{-1})^2}\,.
$$
Therefore
$$
\mu^\infty_\Omega=\lim_{\lambda\to-\infty} \mu_{\lambda}(\Omega)= \lim_{h\to\infty} \mu_{\lambda_h}(\Omega)\ge \mu_0(\mathcal C_\theta)\,.
$$
We proved that $\mu^\infty_\Omega\ge \mu_0(\mathcal C_\theta)$ for any small $\theta>0$, thus $\mu^\infty_0=\infty$
by (\ref{eq:limCtheta}).
\QED

\section{Two inequalities}
\label{A:BV}

Let $w\in H^2_{0,\rm rad}(\ball)$ be a radially symmetric function. Here we prove the following inequalities,
which have been used in Subsection \ref{SS:Auxiliary}:
\begin{gather}
\label{eq:ine_1}
\int_\ball|x|^{-2}|\nabla w|^2\,dx\ge 4\lambda(-\Delta,\disk)\int_\ball |w|^2\,dx\,,\\
\label{eq:ine_2}
\int_\ball|\Delta w|^2\,dx-4\int_\ball|x|^{-2}|\nabla w|^2\,dx\ge \lambda(-\Delta,\disk)\int_\ball |\nabla w|^2\,dx\,.
\end{gather}

\medskip
Thanks to Lemma \ref{L:density} we can assume that $w\in C^\infty_c(\ball\setminus\{0\})$. 

We denote by $r$ and $\rho$ the radii in $\R^4$ and in $\R^2$, respectively. Then, we identify $w\in C^\infty_c(\ball\setminus\{0\})$ with a function $w\in C^\infty_c(0,1)$,
so that $w_r=|x|^{-1}x\cdot \nabla w(x)$ is the radial
derivative of $w$ at $x\in \ball$. Next, we set
$$
g(\rho)=w(\sqrt{\rho})\,.
$$
We  regard at $g\in C^\infty_c(0,1)$ as a radial function $g\in C^\infty_c(\disk\setminus\{0\})$,
so that $g_\rho=|z|^{-1}z\cdot \nabla g(z)$ is the radial derivative of $g$ at $z\in\disk$.
Inequality (\ref{eq:ine_1}) follows from 
$$
\begin{aligned}
\int_\ball|x|^{-2}&|\nabla  w|^2\,dx=|\S^3|\int_0^1 r|w_r|^2\,dr
=2|\S^3|\int_0^1\rho|g_\rho|^2\,d\rho\\
&=\frac{|\S^3|}{\pi} \int_\disk|\nabla g|^2\,dz
\ge \frac{|\S^3|}{\pi}\lambda(-\Delta,\disk)\int_\disk|g|^2\,dz
=4 \lambda(-\Delta,\disk) \int_\ball|w|^2\,dx\,.
\end{aligned}
$$

\medskip
Next, we compute
$$
\begin{aligned}
\frac{1}{|\S^3|}\int_\ball|\Delta w|^2\,dx&=\int_0^1r^3|w_{rr}+3r^{-1}w_r|^2\,dr=
\int_0^1r^3|w_{rr}|^2\,dr+9\int_0^1r|w_{r}|^2\,dr+3\int_0^1r^2(|w_r|^2)_r\,dr\\
&= \int_0^1r^3|w_{rr}|^2\,dr+3\int_0^1r|w_{r}|^2\,dr= \frac{1}{|\S^3|}\Big[\int_\ball|\nabla w_r|^2\,dx+3\int_\ball|x|^{-2}|w_r|^2\,dx\Big]\,.
\end{aligned}
$$
Since $|\nabla w|=|w_r|$, we see that
$$
\int_\ball|\Delta w|^2\,dx-4\int_\ball|x|^{-2}|\nabla w|^2\,dx= \int_\ball|\nabla w_r|^2\,dx-\int_\ball|x|^{-2}|w_r|^2\,dx\,.
$$
Inequality (\ref{eq:ine_2}) follows by applying the improved Hardy inequality (\ref{eq:BV_original}) with $d=4$, $\Omega=\ball$, and \mbox{$w_r\in H^1_0(\ball)$} instead of $u$.
\QED

\medskip
\noindent
{\bf Acknowledgments}.
The second author thanks the Department of Mathematics of the University of Torino 
for the hospitality during her visit in June 2026.


\begin{thebibliography}{XX}
\footnotesize

\bibitem{AH} 
D.~R. Adams and L.~I. Hedberg, 
{\it Function spaces and potential theory}, Grundlehren der mathematischen Wissenschaften, 314, Springer, Berlin, 1996.

\bibitem{Be}
W. Beckner, Weighted inequalities and Stein-Weiss potentials, Forum Math. 20 (2008), 587–606.

\bibitem{BC} 
E.~Berchio and P.~Caldiroli, Weighted Hardy-Rellich inequalities via the Emden-Fowler transform, Preprint
arXiv:2504.09306 (2025).


\bibitem{BM} 
H. Brezis and M.~M. Marcus, Hardy's inequalities revisited, 
Ann. Scuola Norm. Sup. Pisa Cl. Sci. (4) {\bf 25} (1997).


\bibitem{BV} 
H. Brezis and J.~L. V\'azquez, Blow-up solutions of some nonlinear elliptic problems, 
Rev. Mat. Univ. Complut. Madrid {\bf 10} (1997), no.~2, 443--469.

\bibitem{Bor}
D.~I. Borisov and P. Freitas, The spectrum of geodesic balls on spherically symmetric manifolds, 
Comm. Anal. Geom. {\bf 25} (2017), no.~3, 507--544.


\bibitem{C}
C. Cazacu, I. Fidel, Weighted Hardy-Rellich type inequalities: improved best constants and
symmetry breaking, Rev. Roumaine Math. Pures Appl. 69 (2024), 397–413.

\bibitem{CM} 
P. Caldiroli and R. Musina, Rellich inequalities with weights, 
Calc. Var. Partial Differential Equations {\bf 45} (2012), no.~1-2, 147--164.


\bibitem{CM2026} 
P. Caldiroli and R. Musina, 
Existence and nonexistence results for linear biharmonic Dirichlet problems on
unbounded domains with Rellich potential and boundary singularity, preprint (2026).


\bibitem{Fall} 
M.~M. Fall, A note on Hardy's inequalities with boundary singularities, 
Nonlinear Anal. {\bf 75} (2012), no.~2, 951--963.


\bibitem{Fall2} 
M.~M. Fall, On the Hardy-Poincar\'e{} inequality with boundary singularities, 
Commun. Contemp. Math. {\bf 14} (2012), no.~3, 1250019, 13 pp.


\bibitem{FM} 
M.~M. Fall and R. Musina, 
Hardy-Poincar\'e{} inequalities with boundary singularities, 
Proc. Roy. Soc. Edinburgh Sect. A {\bf 142} (2012), no.~4, 769--786.


\bibitem{GGM} 
F. Gazzola, H.-C. Grunau and E.~L. Mitidieri, 
Hardy inequalities with optimal constants and remainder terms, 
Trans. Amer. Math. Soc. {\bf 356} (2004), no.~6, 2149--2168.

\bibitem{GM}
N.~A. Ghoussoub and A. Moradifam, 
Bessel pairs and optimal Hardy and Hardy-Rellich inequalities, 
Math. Ann. {\bf 349} (2011), no.~1, 1--57.



\bibitem{Re}
F. Rellich, Halbbeschr\"ankte Differentialoperatoren h\"oherer Ordnung, 
in {\it Proceedings of the International Congress of Mathematicians, 1954, Amsterdam, vol. III}, pp. 243--250, Erven P. Noordhoff N. V., Groningen.

\bibitem{Ta}
G.~G. Talenti, On the first eigenvalue of the clamped plate, Ann. Mat. Pura Appl. (4) {\bf 129} (1981), 265--280.

\end{thebibliography}
\end{document}